\documentclass{amsart}
\usepackage{amscd}
\usepackage{amssymb}

\input xy
\xyoption{all}
\newcommand{\g}{\mathfrak{g}}
\newcommand{\ba}{\mathfrak{B}}
\newcommand{\gl}{\mathfrak{gl}}
\newcommand{\rip}{\mathfrak{sl}}
\newcommand{\e}{\mathfrak{e}}
\newcommand{\pri}{\mathfrak{c}}
\newcommand{\tri}{\mathfrak{t}}
\newcommand{\K}{\mathbf{K}^\Q}

\newcommand{\iso}{\overset{\sim}{\to}}
\newcommand{\lra}{\longrightarrow}
\newcommand{\liso}{\overset{\sim}{\lra}}
\newcommand{\fib}{\twoheadrightarrow}

\DeclareMathOperator{\tor}{Tor}
\DeclareMathOperator{\pro}{pro}
\DeclareMathOperator*{\holi}{holim}
\DeclareMathOperator*{\coli}{colim}
\DeclareMathOperator{\hofi}{hofiber}
\DeclareMathOperator{\coker}{coker}
\DeclareMathOperator{\tot}{Tot}

\DeclareMathOperator{\primito}{Prim}

\DeclareMathOperator{\ass}{\emph{Ass}}
\DeclareMathOperator{\as1}{\emph{Ass}_1}
\DeclareMathOperator{\rings}{\emph{Rings}}

\newcommand{\E}{{E'}}

\newcommand{\Hy}{\mathbb{H}}
\newcommand{\Q}{\mathbb{Q}}
\newcommand{\Z}{\mathbb{Z}}
\newcommand{\N}{\mathbb{N}}

\newcommand{\C}{\mathcal{C}}

\newcommand{\A}{\mathcal{A}}

\newcommand{\triqui}{\triangleleft}

\theoremstyle{plain}
\newtheorem{thm}{Theorem}[section]
\newtheorem{mth}[thm]{Main theorem}
\newtheorem{cor}[thm]{Corollary}

\newtheorem{lem}[thm]{Lemma}
\newtheorem{prop}[thm]{Proposition}

\theoremstyle{definition}
\newtheorem{defn}[thm]{Definition}
\newtheorem{rem}[thm]{Remark}

\newtheorem{exa}[thm]{Example}

\theoremstyle{remark}
\newtheorem{notation}[thm]{Notation}
\newtheorem*{tute}{Note on notation}
\newtheorem*{ack}{Acknowledgement}
\begin{document}

\title{The
obstruction to excision in $K$-theory and in cyclic homology}

\author[Guillermo Corti\~nas]{Guillermo Corti\~nas*\\{\rm gcorti@dm.uba.ar}}

\address{Departamento de Matem\'atica\\
Ciudad Universitaria Pab 1\\
1428 Buenos Aires\\
Argentina}

\curraddr{Departamento de \'Algebra, Geometr\'\i a y Topolog\'\i a\\
          Facultad de Ciencias, Universidad de Valladolid\\
          Prado de la Magdalena s/n\\
          (47005) Valladolid, Spain.}

\begin{abstract}
Let $f:A\to B$ be a ring homomorphism of not necessarily unital
rings and $I\triangleleft A$ an ideal which is mapped by $f$
isomorphically to an ideal of $B$. The obstruction to excision in
$K$-theory is the failure of the map between relative $K$-groups
$K_*(A:I)\to K_*(B:f(I))$ to be an isomorphism; it is measured by
the birelative groups $K_*(A,B:I)$. Similarly the groups $HN_*(A,B:I)$
measure the obstruction to excision in negative cyclic homology. We show that
the rational Jones-Goodwillie Chern character induces an
isomorphism
$$
ch_*:K_*(A,B:I)\otimes\Q\iso
HN_*(A\otimes\Q,B\otimes\Q:I\otimes\Q).
$$
\end{abstract}

\thanks{(*) Partially supported by CONICET, the ICTP's Associateship and the
Ram\'on y Cajal fellowship and
by grants UBACyT X066, ANPCyT PICT 03-12330 and MTM00958.}

\maketitle
\setcounter{section}{-1}
\section{Introduction}\label{intro}
Algebraic $K$-theory does not satisfy excision. This means that if $f:A\to B$
is a ring homomorphism and $I\triangleleft A$ is an ideal carried
isomorphically to an ideal of $B$, then the map of relative $K$-groups
$K_*(A:I)\to K_*(B:I):=K_*(B:f(I))$ is not an isomorphism in general.
The obstruction
is measured by birelative groups $K_*(A,B:I)$ which are defined so as to fit
in a long
exact sequence
$$
K_{n+1}(B:I)\to K_n(A,B:I)\to K_n(A:I)\to K_n(B:I)
$$
Similarly the obstruction to excision in negative cyclic homology
is measured by birelative groups $HN_*(A,B:I)$. $K$-theory and
negative cyclic homology are related by a character $K_nA\to
HN_nA$, the {\it Jones-Goodwillie Chern character} (\cite[Ch.II]{goo}; see also \cite[\S8.4]{lod}). Tensoring with $\Q$ and
composing with the natural map $HN_n(A)\otimes\Q\to
HN_n(A\otimes\Q)$ we obtain a rational Chern character
\begin{equation}\label{ichee}
ch_n:K_n^\Q(A):=K_n(A)\otimes\Q\to HN_n(A\otimes\Q)=:HN_n^\Q(A)
\end{equation}
The main theorem of this paper is the following.
\begin{mth}\label{mt}
Let $f:A\to B$ be a homomorphism of not necessarily unital rings
and $I\triangleleft A$ an ideal which is carried by $f$
isomorphically onto an ideal of $B$. Then \eqref{ichee} induces an
isomorphism
$$
ch_*:K^\Q_*(A,B:I)\iso HN_{*}^\Q(A,B:I).
$$
If moreover $A$ and $B$ are $\Q$-algebras, then $K_*(A,B:I)$ is a
$\Q$-vectorspace, and the exponent $\Q$ is not needed.
\end{mth}

The theorem above can also be stated in terms of {\it cyclic
homology} (denoted $HC$), as we shall see presently. We recall
that, unlike $HN$, $HC$ commutes with $\otimes\Q$, so that
$HC_*^\Q(A):=HC_*(A\otimes\Q)=HC_*(A)\otimes\Q$.

\begin{cor}\label{xmt}
 There is a natural isomorphism
\[
\nu_*:K_*^\Q(A,B:I)\cong HC^\Q_{*-1}(A,B:I).
\]
\end{cor}
\begin{proof}
Write $HP$ for periodic cyclic homology. There is a long exact
sequence (\cite[5.1.5]{lod})
\begin{equation}\label{nsbi}
HC_{n-1}(A,B:I)\to HN_{n}(A,B:I)\to HP_n(A,B:I) \to
HC_{n-2}(A,B:I).
\end{equation}
Cuntz-Quillen's excision theorem \cite{cq} establishes that
\begin{equation}\label{icq}
HP^\Q_*(A,B:I)=0.
\end{equation}
Here $HP_*^\Q(\ \ ):=HP_*(\ \ \otimes\Q)$. From \eqref{icq} and
from \eqref{nsbi} applied to $A\otimes\Q$, $B\otimes\Q$ and
$I\otimes\Q$, it follows that $HC_{*-1}^\Q(A,B:I)=HN^\Q_*(A,B:I)$.
\end{proof}

Next we shall review some related results in the literature, so as
to put ours in perspective. Bass (\cite[Thm. XII.8.3]{bass}; see
also \cite[pp. 295-298]{kav}) proved $K$-theory satisfies excision
in nonpositive degrees; in our setting this means $K_n(A,B:I)=0$
for $n\le 0$. The analogous result for the negative degrees of
cyclic homology is also true. In fact, by definition
(\cite[2.1.15]{lod}) if $R$ is a unital ring and $J\triqui R$ an ideal then
$HC_{n}(R:J)=0$ for $n<0$ and $HC_0(R,J)=J/[R,J]$, the quotient by
the subgroup generated by the commutators $[r,j]=rj-jr$. Hence in
the situation of \ref{mt}
$$
HC_{-1}(A,B:I)=\coker(I/[A,I]\to I/[B,I])=0.
$$
Thus \ref{xmt} is true for $*$ nonpositive, as both birelative
groups vanish. The particular case of \ref{xmt} when $*=1$ was
proved in \cite{gw0}. In \cite{grw} the statement of Corollary
\ref{xmt} was conjectured to hold when $A$ and $B$ are commutative
unital $\Q$-algebras and $B$ is a finite integral extension of
$A$, ($KABI$ conjecture) and it was shown its validity permits
computation of the $K$-theory of singular curves in terms of their
cyclic homology and of the $K$-theory of nonsingular curves. In
\cite{gw} \ref{xmt} was conjectured for unital $\Q$-algebras and
it was shown that for $*=2$ the left hand side maps surjectively
onto the right hand side.

\smallskip

A special case of the main theorem concerns the birelative groups
$K_*(A:I,J)$ associated to any pair of ideals $I,J\triangleleft
A$; they are defined so as to fit in a long exact sequence
$$
K_{n+1}(A/J:I+J/J)\to K_n(A:I,J)\to K_n(A:I)\to K_n(A/J:I+J/J).
$$
Note that if $I\cap J=0$ then $A\to A/J$ maps $I$ isomorphically
onto $I+J/J=I/I\cap J$, whence \ref{xmt} applies, and we have an a
rational isomorphism

\begin{equation}\label{kij}
K_*^\Q(A:I,J)\iso HC^\Q_{*-1}(A:I,J).
\end{equation}

For $*\le 1$ it is well-known that both birelative $K$- and cyclic
homology groups vanish, whence \eqref{kij} follows. The case $*=2$
was proved independently in \cite{gl} and \cite{keu}. One can
further generalize this to the case when $I\cap J$ is nilpotent,
using a theorem of Goodwillie's (\cite[Main Thm.]{goo} see also
\cite[ \S 11.3]{lod}), which says that if $I\triangleleft R$ is a
nilpotent ideal of a ring $R$, then there is a natural isomorphism
\begin{equation}\label{goowi}
K_*^\Q(R:I)\cong HC^\Q_{*-1}(R:I).
\end{equation}
We point out that, even if Goodwillie states his theorem for $R$
unital, the nonunital case follows from the unital (\cite[\S 4.2]{crelle};
see also Lemma \ref{exnuni} below). As in the case of
\ref{xmt}, \eqref{goowi} can also be stated in terms of negative
cyclic homology. Actually Goodwillie shows (see \cite[0.3]{goo})
that $ch_*$ induces an isomorphism
\begin{equation}\label{goocha}
ch_*:K_*^\Q(R:I)\iso HN^\Q_*(R:I)
\end{equation}
and then uses the singly relative version of \eqref{nsbi} in
combination with another theorem of his, (\cite[II.5.1]{pregoo};
see also \cite[3.5]{cq}), which says that $HP^\Q_*(R,I)=0$ if $I$
is nilpotent. Now we use \eqref{goowi} to generalize \eqref{kij}.

\begin{cor}\label{below}
Let $A$ be a ring and $I,J\triangleleft A$ ideals such that
$I\cap J$ is nilpotent. Then there is a rational isomorphism of birelative
groups
$$
K_*^\Q(A:I,J)\iso HC_{*-1}^\Q(A:I,J).
$$
\end{cor}

\noindent{\it Proof.} The case $I\cap J=0$ is explained above. To
prove the general case, consider the intermediate groups
$K^{inf}_n$ fitting in the long exact sequence
\begin{equation}\label{fib}
\xymatrix{HN_{n+1}^\Q(R)\ar[r]& K^{inf}_{n}(R)\ar[r]&
K_n^\Q(R)\ar[r]^(.45){ch_n}& HN^\Q_n(R)}
\end{equation}
(This notation will be justified below, see \eqref{dkinfi}). By
\ref{mt}, $K^{inf}$ satisfies excision; by \eqref{goocha}, it is
invariant under nilpotent extensions, or {\it nilinvariant}. But a
diagram chase shows that any homology theory of rings $H$
satisfying both excision and nilinvariance verifies $H_*(A:I,J)=0$
if $I\cap J$ is nilpotent. Applying this to both $K^{inf}$ and
$HP^\Q$, we obtain
\begin{equation*}
K_*^\Q(A:I,J)\cong HN_*^\Q(A:I,J)\cong HC^\Q_{*-1}(A:I,J)\qed
\end{equation*}

The result of the corollary above for $A$ unital was announced in
\cite{ow}; however the proof in {\it loc.cit.} turned out to have
a gap (see \cite[pp. 591, line 1]{gw2}). An application of
\ref{below} to the computation of the $K$-theory of particular
rings --other than coordinate rings of curves-- in terms of their
cyclic homology was given in \cite[Thm. 3.1]{gw2};
 see also \cite[Thm. 7.3]{grw}.

As another precedent of the main theorem of this paper, we must cite the work
of Suslin and Wodzicki. To state their theorems we introduce some
notation. We say that a ring $I$ is {\it excisive} for a homology
theory $H$ if $H_*(A,B;I)=0$ for every homomorphism $A\to B$ as in
\ref{mt}. In \cite{wod}, M. Wodzicki characterized those rings
which are excisive for cyclic homology as those whose bar homology
vanishes
\begin{equation}\label{wody}
A \text{ is } HC^\Q-\text{excisive }\iff H^{bar}_*A\otimes\Q=0.
\end{equation}
In fact \eqref{wody} is a particular case of \cite[(3)]{wod}. He
also showed that if $A$ is excisive for rational $K$-theory then
it is excisive for rational cyclic homology (\cite[(4)]{wod}) and
conjectured that the converse also holds. The latter was proved by
Suslin and Wodzicki; they showed (\cite[Thm. A]{qs})

\begin{equation}\label{tqs}
H^{bar}_*(A)\otimes\Q=0\Rightarrow A \text{ is excisive for
$K^\Q$-theory.}
\end{equation}

They proved further that if $A$ is a $\Q$-algebra then \eqref{tqs} still
holds even if we do not tensor with $\Q$ (\cite[Thm. B]{qs}). Note that \eqref{tqs} is a
formal consequence of \ref{xmt} and \eqref{wody}. Actually our
proof of \ref{mt} involves proving a version of Suslin-Wodzicki's
theorem for a certain type of pro-rings (Theorem \ref{gqs}).

\smallskip

\noindent{\it Sketch of Proof of \ref{mt}}. The assertion
concerning $\Q$-algebras follows (see \ref{chum} below) from
Weibel's result \cite{chu} (see also \cite[1.9]{qs}) that for all
$n\ge 2$, $K$-theory with $\Z/n$ coefficients satisfies excision
for such algebras. The rest of the proof has four parts:

\smallskip

\noindent{a)} {\it An abstraction of arguments of Cuntz and
Quillen  and its combination with the tautological characters of
\cite{crelle}}. (This is done in Section \ref{cqp} below). We consider functors
from the category $\ass_1=\ass_1(k)$ of unital rings
over a commutative ring $k$ to fibrant spectra which preserve products up to homotopy.
If $X$ is such a functor and $I\triqui A\in\ass_1$ an ideal, we put
\begin{gather*}
\ X(A:I):=\hofi(X(A)\to X(A/I)),\ \ \hat{X}(A/I^\infty):=\holi_n
X(A/I^n),\\ \hat{X}(A:I^\infty):=\holi_n X(A:I^n).
\end{gather*}

If $f:A\to B\in \ass_1$ is as in \ref{mt}, we set
\begin{align*}
X(A,B:I)=&\hofi(X(A:I)\to X(B:I)),\\
\hat{X}(A,B:I^\infty)=&\hofi(\hat{X}(A:I^\infty)\to\hat{X}(B:I^\infty)).
\end{align*}
We call $I$ $\infty$-excisive if $\hat{X}(A,B:I^\infty)$ is weakly
contractible for every homomorphism $A\to B$ as above. We say that
$X$ is {\it nilinvariant} if $X(A:I)$ is contractible for every
nilpotent ideal $I\triqui A$. We show that if $X$ is nilinvariant
and every ideal $I$ of every free unital algebra is
$\infty$-excisive then $X$ is excisive (Proposition \ref{exci}).
This gives a criterion for proving excision which generalizes that
used by Cuntz and Quillen in the particular case of $HP$ in
\cite{cq}. Next, given any not necessarily excisive or
nilinvariant functor $X$ as above, we consider its {\it
noncommutative infinitesimal hypercohomology} (\cite[\S5]{crelle})
$\Hy(A_{inf},X)$. This is a nilinvariant functor and is equipped
with a natural map $\Hy(A_{inf},X)\to X(A)$. Write $\tau X(A)$ for
the delooping of the fiber of the latter map. We have a homotopy
fibration
\begin{equation}\label{hinfi}
\Hy(A_{inf},X)\to X(A)\overset{c^\tau}\lra \tau X(A).
\end{equation}
We call $c^\tau$ the {\it tautological character}. Applying the
Cuntz-Quillen criterion to $\Hy(A_{inf},X)$, we obtain (Theorem
\ref{fund}) that if every ideal of every free unital algebra
$F\in\ass_1$ is $\infty$-excisive for both $X$ and $\tau X$ then
the map $X(A,B:I)\to \tau X(A,B:I)$ is a weak equivalence for
every algebra homomorphism as in \ref{mt}.
\smallskip

\noindent{b)}{\it Agreement of the tautological and rational
Jones-Goodwillie characters.} (Section \ref{checa}). We apply
\eqref{hinfi} to $X=\K$, the nonconnective rational $K$-theory
spectrum. We show that there is a natural isomorphism
$HN_n^\Q(A)\cong \tau K^\Q_n(A)$ under which
the tautological character $c^\tau_n$ is identified with $ch_n$
(Theorem \ref{agree}). In particular, the {\it infintesimal $K$-theory groups}
\begin{equation}\label{dkinfi}
K^{inf}_nA:=\pi_n\Hy (A_{inf},\K) \qquad (n\in\Z).
\end{equation}
are the ``intermediate groups" of \eqref{fib}.
Theorem \ref{agree} is of independent interest as it
shows that $HN^\Q$ can be functorially derived from $K^\Q$,
and extends to all $n\in \Z$ previous results of the author
(\cite[6.3]{chino}; \cite[6.2]{crelle}; \cite[5.1]{sheaf}).

\smallskip

\noindent{c)}{\it Suslin-Wodzicki's theorem for
$A^\infty$}. (Section \ref{sgqs}). We show that if $A$ is a ring
such that for all $r\ge 0$ the pro-vectorspace
$$
H^{bar}_r(A^\infty)\otimes\Q:=\{H^{bar}_r(A^n)\otimes\Q\}_n
$$
is zero, then $A$ is $\infty$-excisive for $K^\Q$ (Theorem
\ref{gqs}). Note that in the particular case when $A=A^2$, the
latter assertion and \eqref{tqs} coincide, since $A^\infty$ is
just the constant pro-ring $A$ in this case. Our proof follows the
strategy of Suslin and Wodzicki's proof of \eqref{tqs}, (see the
summary at the beginning of Section \ref{sgqs}) and adapts it to
the pro-setting. Some technical results on pro-spaces
needed in this section are proved in the
Appendix.

\smallskip

\noindent{d)} {\it Application of known results on bar and cyclic
pro-homology}. (Section \ref{qsp}). By parts (a), (b) and (c), to
finish the proof it is enough to show that if $I$ is an ideal of a
free unital ring, then both (i) and (ii) below hold.

\smallskip

\noindent{(i)} $H^{bar}_r(I^\infty)\otimes\Q=0$ ($r\ge 0$).

\smallskip

\noindent{(ii)} $I$ is $\infty$-excisive for $HN^\Q$.

\smallskip

By \eqref{nsbi} and \eqref{icq} the latter property is equivalent
to

\smallskip

\noindent{(ii)'} $I$ is $\infty$-excisive for $HC^\Q$.

\smallskip

As $\otimes\Q$ commutes with both $H^{bar}$ and $HC$ and sends
free unital rings to free unital $\Q$-algebras, it suffices to
verify (i) and (ii)' for ideals of free unital $\Q$-algebras. Both
of these are well-known and are straightforward from results in
the literature; see 4.2 for details.\qed

\bigskip

The rest of this paper is organized as follows. In Section
\ref{cqp} we carry out part (a) of the sketch above. We generalize
Cuntz-Quillen's excision principle (Proposition \ref{exci}),
recall the construction of the tautological character (\ref{infx})
and obtain a criterion for proving that the latter computes the
obstruction to excision (Theorem \ref{fund}). Part (b) corresponds
to Section \ref{checa}, where we show that the tautological
character for rational $K$-theory is the rational Jones-Goodwillie
character (Theorem \ref{agree}). Section \ref{sgqs} is devoted to
part (c); we prove a version of \eqref{tqs} for pro-rings of the
form $A^\infty$ (Theorem \ref{gqs}). The proof of \ref{mt} is
completed in Section \ref{qsp}, where we carry out part (d) of the
sketch and show that if in the situation of the Main Theorem, $A$
and $B$ are $\Q$-algebras, then the groups $K_*(A,B:I)$ are
$\Q$-vectorspaces (Lemma \ref{chum}). Notation for pro-spaces as
well as some technical results on them which are used in Section
\ref{sgqs} are the subject of the Appendix.

\begin{tute} In this paper space=simplicial set.
\end{tute}

\begin{ack} I learned about the existence of the $KABI$ conjecture in
1991
from C. Weibel. His papers with Geller and Reid (\cite{grw}) and
Geller (\cite{gw0}, \cite{gw},\cite{gw2}) demonstrated its
potential applications, produced first positive results and thus
made it attractive as a problem. I am thankful to all three of
them for calling my attention to it through their work, as well as
to J. Cuntz, D. Quillen, A. Suslin and M. Wodzicki, for their
results are crucial to this paper. Special thanks to the editor
and the referees for the work and dedication spent in improving
my paper.
\end{ack}

\section{Cuntz-Quillen excision principle}\label{cqp}
\noindent{\sc Summary.} In this section we formulate in an
abstract setting the method used by Cuntz and Quillen to prove
excision for periodic cyclic homology of algebras over a field of
characteristic zero. In their proof they first show that $HP$
satisfies excision for pro-ideals of the form $I^\infty$ where $I$
is an ideal of a free algebra and then combine this with the
invariance of $HP$ under nilpotent extensions to prove excision
holds for all ideals. Our setting is that of functors $X$ from the
category of algebras over a commutative ring $k$ to that of
fibrant spectra. We assume that $X$ preserves products up to homotopy.
We show that if $X$ is invariant (up to weak
equivalence) under nilpotent extensions and satisfies excision for
pro-rings of the form $I^\infty$ whenever $I$ is an ideal of a
free unital algebra, then $X$ satisfies excision for all algebras
(Prop. \ref{exci}). Then we apply this to give a criterion for
proving that the tautological character $c^\tau:X\to \tau X$ of
\cite{crelle} induces an equivalence at the level of the groups
which measure the obstruction to excision (Theorem \ref{fund}).
$\qed$

\bigskip

We consider associative, not necessarily unital algebras over a
fixed ground ring $k$. We write $\ass:=\ass(k)$ for the category
of algebras, $\as1$ for the subcategory of unital algebras and
unit preserving maps, and
\begin{equation}\label{tilde}
A\mapsto \tilde{A}
\end{equation}
for the left adjoint functor to the inclusion $\as1\subset\ass$.
By definition,
$$
\tilde{A}=A\oplus k,\ \ (a,\lambda)(b,\mu):=(ab+\mu a+\lambda
b,\lambda\mu).
$$
Throughout this section $X$ will be a fixed functor from $\as1$ to
fibrant spectra (terminology for spectra is as in \cite{tho}). We shall assume
that $X$ preserves finite products up to homotopy; this means that if $A,B\in\as1$ then
the canonical map is a weak equivalence:
\begin{equation}\label{prodpres}
X(A\times B)\iso X(A)\times X(B).
\end{equation}

Let
$A\in\as1$, $I\triangleleft A$ an ideal, $\pi:A\to A/I$ the
projection. Write $X(A:I)$ for the homotopy fiber of $X(\pi)$. We
say that $A\in\ass$ is {\it excisive} (with respect to our fixed functor
$X$) if for every
commutative diagram
\begin{equation}\label{diag}
\xymatrix{&A\ar@{_{(}->}[dl]\ar@{^{(}->}[dr]&\\
          B\ar[rr]&&C}
\end{equation}
such that $B\to C$ is a map in $\as1$ and both arrows out of $A$ embed it as
an ideal, the map
$$
X(B:A)\lra X(C:A)
$$
is an equivalence. Equivalently, $A$ is excisive if
$$
X(B,C:A):=\hofi(X(B:A)\to X(C:A))
$$
is contractible for every diagram \eqref{diag}. We note it is
enough to check this condition for $B=\tilde{A}$. The functor $X$
is called {\it excisive} if every $A\in\ass$ is.

The rule
\begin{equation}\label{ye}
A\to Y(A):=X(\tilde{A}:A)
\end{equation}
defines a functorial fibrant spectrum on $\ass$. If $A\in\as1$,
then $\tilde{A}\cong A\times k$, and because the homotopy category
of spectra is additive, we have
$$
X(\tilde{A})\iso Y(A)\times X(k).
$$
Since we are assuming that $X$ preserves finite products up to homotopy \eqref{prodpres}, it follows
that for $A\in \as1$, there is a natural homotopy equivalence
$$
Y(A)\iso X(A).
$$
Thus $Y$ extends $X$ up to homotopy. Note further
that, by definition,
$$
Y(0)=\hofi(X(k)\overset{1}\lra X(k))
$$
is contractible. If $A\triqui B$ is an ideal we put
$$
Y(B:A):=\hofi(Y(B)\to Y(B/A)).
$$
The natural map $Y(R)\to X(\tilde{R})$ induces
\begin{equation}\label{ytx}
Y(B:A)\to X(\tilde{B}:A).
\end{equation}
\begin{lem}\label{exnuni}
The map \eqref{ytx} is a weak equivalence. In particular,
$Y(A:A)\iso YA$.
\end{lem}
\begin{proof}
Consider the commutative diagram
$$
\xymatrix{Y(B:A)\ar[d]\ar[r]&YB\ar[d]\ar[r]&Y(B/A)\ar[d]\\
          X(\tilde{B}:A)\ar[d]\ar[r]&X(\tilde{B})\ar[d]\ar[r]&X(\widetilde{B/A})\ar[d]\\
          {*}\ar[r]&X(k)\ar[r]^1&X(k)}
$$
Here $1$ denotes the identity map. By definition, all rows as well
as the middle and right hand columns are homotopy fibrations; it
follows that also the column on the left is one.
\end{proof}
\begin{cor}\label{exnuni1}
Let $B\to C$ be a homomorphism in $\ass$ carrying an ideal
$A\triangleleft B$ isomorphically onto an ideal of $C$. Then the
natural map
$$
Y(B,C:A):=\hofi(Y(B:A)\lra Y(C:A))\to X(\tilde{B},\tilde{C}:A)
$$
is an equivalence.
\end{cor}
\begin{cor}\label{exnuni2}
If $A$ is excisive and $A\triqui B$, then the map $Y(A)\to Y(B:A)$
is a weak equivalence.
\end{cor}
\begin{proof}
By \ref{exnuni}, it suffices to show that
\begin{equation}\label{ese}
Y(A:A)\to Y(B:A)
\end{equation}
is a weak equivalence. But by \ref{exnuni1}, the homotopy fiber of
\eqref{ese} is weakly equivalent to $X(\tilde{A},\tilde{B}:A)$,
which is contractible because $A$ is excisive.
\end{proof}

\begin{lem}\label{ifree}
If every ideal of every free unital algebra is excisive, then $X$
is excisive.
\end{lem}
\begin{proof} Let $F:\emph{Sets}\to\as1$ be the free unital algebra functor.
For $R\in\as1$ let $JR=\ker(FR\to R)$ and $\rho:R\to FR$ be
respectively the kernel of the projection and its natural
set-theoretic section. The ideal $JR$ is generated by $\rho(1)-1$,
by $\rho(0)-0$ and by the elements of the form
$\rho(a+b)-\rho(a)-\rho(b)$ and of the form
$\rho(ab)-\rho(a)\rho(b)$ ($a,b\in R$). In particular $J$
preserves surjections. It follows from this and the snake lemma
that if $B\in\as1$ and $A\triangleleft B$ is an ideal, then
$I=\ker(F(B)\to F(B/A))$ maps onto $A$. The kernel of the
surjection $\tilde{I}\to \tilde{A}$ is $L:=JB\cap I$. We have a
commutative diagram
\begin{equation}\label{grandiagrama}
\xymatrix@!0{&JB\ar[rr]\ar'[d][dd]&&J(B/A)\ar[dd]\\
L\ar[ur]\ar[rr]\ar[dd]&&0\ar[ur]\ar[dd]\\
&FB\ar[rr]\ar'[d][dd] && F(B/A)\ar[dd]\\
\tilde{I}\ar[ur]\ar[rr]\ar[dd]&&k\ar[ur]\ar[dd]\\
&B\ar'[r][rr] &&B/A\\
\tilde{A}\ar[rr]\ar[ur]&& k\ar[ur]}
\end{equation}
Let $Y$ be as in \eqref{ye}. By \ref{exnuni2}, applying $Y$ to the
diagram above yields a diagram whose columns are homotopy
fibrations, and of which the two top squares are homotopy
cartesian. It follows that also the bottom square is homotopy
cartesian.
\end{proof}

We say that $X$ is {\it nilinvariant} if $X(B:A)$ is contractible
for every nilpotent ideal $A\triangleleft B$. If $A=\{A_n\}$ is a
pro- unital algebra, and $I=\{I_n\}$ a (levelwise) pro-ideal we
put
$$
\hat{X}(A)=\holi_n X(A_n),\quad \hat{X}(A:I)=\holi_n X(A_n:I_n).
$$
We say that $A$ is
$\infty$-excisive if for every diagram \eqref{diag} the induced
map
$$
\hat{X}(B:A^\infty)\lra \hat{X}(C:A^\infty)
$$
is an equivalence. Here as in \cite{cq}, $A^{\infty}=\{A^n\}_n$
with the natural inclusions as transition maps.

\begin{lem}\label{nilex}
If $X$ is nilinvariant and $A$ is $\infty$-excisive, then $A$ is excisive.
\end{lem}
\begin{proof}
Let $B\to C$ be as in \eqref{diag}. We have a commutative diagram
of level pro-maps
$$
\begin{CD}
B@>>>B/A^\infty @>>>B/A\\
@VVV @VVV @VVV\\
C @>>>C/A^\infty @>>>C/A\\
\end{CD}
$$
Apply $\hat{X}$ to this diagram, noting that $\hat{X}(R)\iso X(R)$ if
$R$ is a constant pro-algebra.
\end{proof}

\begin{prop}\label{exci}
Assume that $X$ is nilinvariant and that every ideal of every free
unital algebra is $\infty$-excisive. Then $X$ is excisive.
\end{prop}
\begin{proof}
Immediate from \ref{nilex} and \ref{ifree}.
\end{proof}

\begin{defn}\label{infx}
The {\it infinitesimal hypercohomology} of an object $A\in\as1$ with
coefficients in $X$ is
the fibrant spectrum
\begin{equation}\label{finfx}
\Hy(A_{inf},X):=\holi_{A_{inf}}X.
\end{equation}
Here $A_{inf}$ is the category of all surjections $B\fib A\in\ass$
with nilpotent kernel; $X$ is viewed as a functor on $A_{inf}$ by
$(B\to A)\mapsto X(B)$. Although the category $A_{inf}$ is large,
it is proved in \cite[5.1]{crelle} that it has a left cofinal small
subcategory in the sense of \cite[Ch IX\S9]{bk}, which by
\cite[2.2.1]{crelle} can be chosen to depend functorially on $A$.
The homotopy limit is taken over this small subcategory.
  There is a natural map $\mathbb{H}(A_{inf},X)\to X(A)$;
we write $\tau X$ for the delooping of its homotopy fiber.
We have a map of spectra
$$
c^\tau:X(A)\to \tau X (A)
$$
of which the homotopy fiber is weakly equivalent to
$\mathbb{H}(A_{inf},X)$. We call $c^\tau$ the {\it tautological
character} of $X$.
\end{defn}
\begin{thm}\label{fund}
Assume every ideal of every free unital algebra is
$\infty$-excisive for both $X$ and $\tau X$. Let $A\to B\in\as1$
be a homomorphism carrying an ideal $I\triangleleft A$
isomorphically onto an ideal of $B$. Then the tautological
character induces a weak equivalence
$$
c^\tau:X(A,B:I)\iso \tau X(A,B:I).
$$
In particular the functor $A\mapsto \Hy(A_{inf},X)$ is excisive.
\end{thm}

\begin{proof} Immediate from \ref{exci} and the nilinvariance of
$\mathbb{H}(A_{inf},X)$ (see \cite[5.2]{crelle}).
\end{proof}
\begin{rem}\label{exnuni3}
For $A\in\ass$, the tautological character $c^\tau:X\to \tau X$
induces a map
\begin{equation}\label{tau'}
c^{\tau'}:Y(A)\to \tau'Y(A):=\tau X(\tilde{A}:A).
\end{equation}
If $A$ happens to be unital, we have a commutative diagram with
vertical weak equivalences
\begin{equation}\label{lindo}
\xymatrix{X(A)\ar[d]_{\wr}\ar[r]^{c^{\tau'}}& \tau X(A)\ar[d]^{\wr}\\
Y(A)\ar[r]_{c^{\tau'}}& \tau'Y(A).}
\end{equation}
Thus $c^\tau$ and $c^{\tau'}$ agree up to homotopy for $A\in\as1$.
By \ref{exnuni1}, under the hypothesis of \ref{fund},
$$
c^{\tau'}:Y(A,B:I)\to \tau' Y(A,B:I)
$$
is an equivalence. On the other hand, if in
formula \eqref{finfx} we substitute for $A_{inf}$ the category
$inf(\ass\downarrow A)$ of all nilpotent extensions of $A$ in
$\ass$, and define $\tau YA$ as the delooping of the homotopy
fiber of the canonical map, we obtain a homotopy fibration
sequence
\begin{equation}\label{tauy}
\Hy(inf(\ass\downarrow A),Y)\to Y(A)\overset{c^\tau}\lra \tau Y(A).
\end{equation}
The following Lemma says that the characters of \eqref{tauy} and
\eqref{lindo} agree up to homotopy; this will be used in Section
\ref{checa} (Proof of \ref{agree} and Remarks \ref{chetaununi} and
\ref{chetauq}).

\end{rem}

\begin{lem}\label{cotau}
Let $A\in\ass$. Then there is a weak equivalence $\tau'Y(A)\iso \tau Y(A)$
such that the diagram
$$
\xymatrix{&Y(A)\ar[dl]_{c^{\tau'}}\ar[dr]^{c^{\tau}}& \\
         \tau'Y(A)\ar[rr]^\sim&&\tau Y(A)}
$$
is homotopy commutative.
\end{lem}
\begin{proof}
We shall abuse notation and write $A_{inf}$ for
$inf(\ass\downarrow A)$ and $\tilde{A}_{inf}$ for
$inf(\as1\downarrow\tilde{A})$. Let
$FA:=\hofi(\Hy(\tilde{A}_{inf},X)\to \Hy(k_{inf},X))$; there is a
homotopy fibration sequence
\[
FA\to Y(A)\to \tau'Y(A).
\]
The functor
\begin{equation}\label{tilinfi}
\widetilde{\ \ }:A_{inf}\to \tilde{A}_{inf}
\end{equation}
is left adjoint to the pullback over $A\to \tilde{A}$. Hence \eqref{tilinfi}
is left cofinal in the sense of \cite{bk}, whence
$$
\Hy(\tilde{A}_{inf},X)\to\Hy( A_{inf},X\circ\widetilde{\ \ })
$$
is a weak equivalence. Thus in the following diagram with homotopy fibrations as rows,
$$
\xymatrix{FA\ar[d]\ar[r]&\Hy(\tilde{A}_{inf},X)\ar[d]\ar[r]& \Hy(k_{inf},X)\ar[d]\\
F'A\ar[r]&\Hy(A_{inf},X\circ\widetilde{\ \ })\ar[r]&
\Hy(0_{inf},X\circ\widetilde{\ \ })}
$$
the columns are weak equivalences. On the other hand we have a commutative diagram

$$
\xymatrix{\Hy(A_{inf},Y)\ar[r]\ar[d]&\Hy(A_{inf},X\circ\widetilde{\ \ })\ar[d]\ar[r]&\Hy(A_{inf},X(k))\ar[d]\\
          \Hy(0_{inf},Y)\ar[r]&\Hy(0_{inf},X\circ\widetilde{\ \ })\ar[r]&\Hy(0_{inf},X(k))}
$$
Because $\Hy(A_{inf},\ \ )=\holi_{A_{inf}}$ preserves homotopy
fibrations, both rows above are homotopy fibrations. Note $F'A$ is
the homotopy fiber of the vertical map in the middle. Let $F"A$ be
the homotopy fiber of the vertical map on the left; I claim that
both $F"A\to \Hy(A_{inf},Y)$ and $F"A\to F'A$ are weak
equivalences. In fact because $0_{inf}$ has an initial object, the
leftmost spectrum in the bottom row is weak equivalent to $Y(0)$,
and therefore contractible, while the other two spectra in the
same row are weak equivalent to each other and to $X(k)$.
To prove that the rightmost vertical arrow is a weak equivalence,
it suffices to show that the natural map
$\Hy(A_{inf},X(k))\to X(k)$ is one. But this is clear from the fact that there is a left cofinal
functor $\N^{op}\times\Delta\to A_{inf}$ (\cite[5.1]{crelle}). The claim is proved. We have
thus constructed a diagram
\begin{equation*}
\xymatrix{FA\ar[r]^{\sim}\ar[drr]&F'A&&F"A\ar[ll]_\sim\ar[r]^(0.4)\sim&\Hy(A_{inf},Y)\ar[dll]\\
                                    &&Y(A)&&}
\end{equation*}
in which all horizontal arrows are weak equivalences. Inverting
$F"A\to F'A$ up to homotopy, taking homotopy fibers of the slanted
arrows and delooping them, we get the homotopy commutative diagram
of the lemma.
\end{proof}

\section{The Chern character is the tautological character}\label{checa}

\noindent{\sc Summary.} In this section we recall the construction
of the rational Chern character of Jones and Goodwillie
\begin{equation}\label{chee}
ch_n:K_n^\Q(A)\to HN^\Q_n(A)
\end{equation}
for nonnegative $n$ and extend it to all $n\in\Z$. On the other
hand the general construction of the previous section applied to
the nonconnective rational $K$-theory spectrum $\K(A)$ gives a
tautological character $c^\tau_n:K_n^\Q(A)\to \tau K^\Q_n(A)$. We
prove that the characters $ch_n$ and $c^\tau_n$ agree up to a
canonical isomorphism $HN_n^\Q(A)\cong\tau K^\Q_n(A)$ $(n\in\Z)$
(see \ref{agree}). In particular this shows that \eqref{chee}
comes from a map of spectra which is functorially derived from
$K^\Q$.$\qed$

\bigskip

Let $A$ be a unital ring. We write $BG$ for the nerve of the group
$G$; thus for us $BG$ is a pointed simplicial set. The rational
plus construction of the general linear group is the Bousfield-Kan
$\Q$-completion
$$
K^\Q(A):=\Q_{\infty}BGL(A).
$$
There is a nonconnective spectrum $\K A$ of which the $n$-th space
is
\begin{equation}\label{nspace}
{}_n\K A:=\Omega K^\Q(\Sigma^{n+1}A)
\end{equation}
where $\Sigma$ is Karoubi's suspension functor (\cite{gersten}).
Next we consider the complex for negative cyclic homology, which
we call $CN_*$; this is the normalized version of that denoted
$\tot_*\ba C^{-}$ in \cite[5.1.7]{lod}; normalization is as in
\cite[2.1.9]{lod}. For $n\ge 1$ the character \eqref{chee} is
induced by a map of spaces
\begin{equation}\label{chespace}
K^\Q (A)\overset{ch}{\to} SCN^\Q_{\ge 1}(A).
\end{equation}
Here $SCN_{\ge 1}$ is the simplicial abelian group the Dold-Kan
correspondence associates to the truncation of $CN_*$ which is $0$
in degree $0$, the kernel of the boundary operator in degree one,
and $CN_n$ in degrees $n\ge 2$. The superscript $\Q$ means
$SCN_{\ge 1}$ is applied to $A\otimes\Q$. The map $ch_0$ can be
obtained from $ch_1$, using the inclusion $A\subset A[t,t^{-1}]$
and naturality; see \cite[8.4.10]{lod}. Consider the sequence
\begin{equation}\label{suse}
0\to M_{\infty}A\to \Gamma A\to \Sigma A\to 0
\end{equation}
where $\Gamma$ is Karoubi's cone functor. As $A$ is unital,
$M_{\infty}A$ is excisive for both $K^\Q$-theory and $HN^\Q$ (cf.
\cite{wod}) and has the same $K^\Q$- and $HN^\Q$-groups as $A$.
Moreover $K^\Q_*(\Gamma A)=HN^\Q_*(\Gamma A)=0$. Thus by
naturality we get a commutative diagram with vertical isomorphisms
\begin{equation}\label{conmuta}
\xymatrix{K_{n}^\Q(A)\ar[rr]^{ch_n} &&HN^\Q_n(A)\\
          K_{n+1}^\Q(\Sigma A)\ar[rr]_{ch_{n+1}\Sigma}\ar[u]^{\cong}
          &&
HN_{n+1}^\Q(\Sigma A)\ar[u]^{\cong}}
\end{equation}
This says that the map
\begin{equation}\label{looch}
ch_{n}:=ch_{n+r}\Sigma^r:K_n^\Q(A)\to HN^\Q_{n}(A)
\end{equation}
does not depend on $r\ge 1-n$. On the other hand the general
framework of the previous section produces a map of spectra
$$
c^\tau:\K(A)\to\mathbf{\tau K}^\Q(A)
$$
which at the level of homotopy groups gives a character
$$
c^\tau_n:K_n^\Q (A)\to \tau K^\Q_n(A)\qquad (n\in \Z).
$$
We shall show that this map agrees with the rational
Jones-Goodwillie Chern character up to canonical isomorphism.

\begin{thm}\label{agree}
There is a natural isomorphism $HN^\Q_n(A)\cong\tau K_n^\Q(A)$
($n\in \Z$) which makes the following diagram commute
$$
\xymatrix{K^\Q_n(A)\ar[r]^{c^\tau_n}\ar[dr]^{ch_n} &\tau K_n^\Q(A)\\
           &HN_n^\Q(A)\ar[u]^{\wr}}
$$
\end{thm}
\begin{proof}
By \cite[5.1 and 5.5]{sheaf}, for every $R\in \rings:=\ass(\Z)$
there is a homotopy fibration sequence of spaces
\begin{equation}\label{fisi}
\xymatrix{ \mathbb{H}(R_{inf},K^\Q)\ar[r]
&K^\Q(R)\ar[r]^(0.4){ch}&
 SCN_{\ge 1}^\Q(R).}
\end{equation}
Here $R_{inf}=inf(\rings\downarrow R)$ is the nonunital version of
the infinitesimal site, and infinitesimal hypercohmology is taken with coefficients in a space,
rather than a spectrum; this is defined by the same formula \eqref{finfx} (see \cite[\S5]{crelle}).
By definition, $K^\Q R:=\hofi (K^\Q\tilde{R}\to K^\Q\Z)$.
Note that $CN_{\ge 1}^\Q\Z=CN_{\ge 1}\Q=0$, so that $SCN^\Q_{\ge 1}R:=\hofi (SCN^\Q_{\ge 1}\tilde{R}\to SCN^\Q_{\ge
1}\Z)\cong SCN^\Q_{\ge 1}\tilde{R}$. Using \eqref{fisi} we will
prove that for $A$ unital, $ch_n$ agrees with the map induced by
the nonunital version \eqref{tauy} of $c^\tau$; by \eqref{lindo}
and Lemma \ref{cotau}, the latter map is the same as that induced
by the unital version. Applying \eqref{fisi} to $R=\Sigma^{n+1}A$
we get the homotopy fibration sequence in the right column of the
following diagram

\begin{equation}\label{25}
\xymatrix{\Hy (A_{inf},\Hy (\Sigma^{n+1}(\ \
)_{inf},K^\Q))\ar[d]\ar[r]&
\mathbb{H}(\Sigma^{n+1}A_{inf},K^\Q)\ar[d]\\
\mathbb{H}(A_{inf},K^\Q\Sigma^{n+1})\ar[r]\ar[d]&K^\Q(\Sigma^{n+1}A)\ar[d]^(0.4){ch\Sigma^{n+1}}\\
\mathbb{H}(A_{inf}, SCN_{\ge 1}^\Q\Sigma^{n+1})\ar[r]&
SCN_{\ge 1}^\Q(\Sigma^{n+1}A)}
\end{equation}

The column on the left results from applying
$\mathbb{H}(A_{inf},\ \ )$ to the homotopy fibration sequence on the
right; because $\holi$ preserves such sequences, it is again one.
It follows from \cite[\S5]{crelle} and from the fact that $\Sigma$
preserves nilpotent  extensions that the horizontal map at the top
is an equivalence. Thus the map between loopspaces of the fibers
of the first two horizontal maps from the bottom is an
equivalence.
 Because loopspaces commute with homotopy limits, and in view of the identity
\eqref{nspace}, we get that the vertical map to the left of the following diagram is an equivalence
$$
\xymatrix{{}_{n-1}\tau \K
(A)\ar[r]\ar[d]^{\wr}&\Hy(A_{inf},{}_{n}\K
)\ar[r]\ar[d]&{}_n\K(A)\ar[d]\\
F \ar[r]&\Hy(A_{inf},\Omega SCN_{\ge 1}^\Q\Sigma^{n+1})\ar[r] &\Omega SCN_{\ge 1}^\Q(\Sigma^{n+1}A)}
$$
Here both rows are homotopy fibrations and ${}_{n-1}\tau \K$ is
the $n-1$ space of the spectrum $\tau\K$.  We shall show that
$\mathbb{H}(A_{inf},\Omega SCN_{\ge 1}^\Q\Sigma^{n+1})$ is
contractible; the lemma follows from this. Because $HN$ is Morita
invariant and satisfies excision for sequences of the form
\eqref{suse}, there is an equivalence
\begin{equation}\label{snsig=cn-}
\mathbb{H}(A_{inf},\Omega SCN_{\ge 1}^\Q(\Sigma^{n+1}(\ \
)))\overset{\sim}{\to} \mathbb{H}(A_{inf}, SCN_{\ge -n}^\Q[-n]).
\end{equation}
Here $S$ is again the Dold-Kan correspondence and $CN_{\ge
-n}[-n]$ is the negative cyclic complex which we have truncated
below $n$ in the same way as explained above for $n=1$, and then
shifted upwards to a nonnegative chain complex. By \cite[5.32]{tho}
it suffices to show that the infinitesimal hypercohomology of the
truncated complex $CN_{\ge -n}[-n]$ vanishes in nonnegative
degrees. By \cite[5.4]{sheaf}, we can replace $CN_*$ by the complex
$NX_*$ for the negative cyclic homology of Cuntz-Quillen's mixed
complex $X_*$. We have a second quadrant spectral sequence
$$
E^1_{p,q}=H^{-p}(A_{inf},(NX^\Q_{\ge -n})_{q-n})
\Rightarrow \Hy^{-p-q}(A_{inf},NX^\Q_{\ge -n}[-n])
$$
By the proof of \cite[5.1]{sheaf}, $E^1_{p,q}=0$ for $q\neq 0$ and also for $q=0$ and $p\ge -2$.
\end{proof}

\begin{rem}\label{chetaununi}
If $A$ is a nonunital ring, the Chern character is defined as the
map
\begin{equation}\label{nuch}
K^\Q_*(A):=K^\Q_*(\tilde{A}:A)\to HN^\Q_*(\tilde{A}:A)=:HN^\Q_*(A)
\end{equation}
induced by $ch_*:K^\Q_*(\tilde{A})\to HN^\Q_*(\tilde{A})$. By the
theorem above, \eqref{lindo}, and Lemma \ref{cotau}, \eqref{nuch}
coincides with the map induced by the tautological character
\eqref{tauy}.
\end{rem}
\begin{rem}\label{chetauq}
If $A\in\ass(\Q)$, by \cite[5.1]{sheaf}, one has a homotopy
fibration sequence
\begin{equation}\label{intsheaf}
\xymatrix{ \mathbb{H}(A_{inf},K)\ar[r] &K(A)\ar[r]^(0.4){ch^\Z}&
 SCN_{\ge 1}(A)}
\end{equation}
where $K=\Z_{\infty}GL$ is the integral plus construction and
$ch^\Z$ is the integral Jones-Goodwillie character. By the same
method as described above in this section for $ch_*=ch^\Q_*$, one
can define characters $ch^\Z_*$ for all $*\in\Z$. If further $A$
is unital, using \eqref{intsheaf} and the same argument as that of
the proof of \ref{agree}, one gets that $ch^\Z_*$ agrees with the
tautological character for the integral nonconnective $K$-theory
spectrum $\mathbf{K}$, obtained by delooping the fiber of
$\Hy(inf(\ass(\Q)\downarrow A),\mathbf{K})\to \mathbf{K}$, as well
as with that obtained using $inf(\as1(\Q)\downarrow A)$ instead.
This can be further extended to nonunital $A$ by the same
argument as in \ref{chetaununi}.
\end{rem}

\section{Suslin-Wodzicki theorem for $A^\infty$}\label{sgqs}

\noindent{\sc Summary.} In this section we show (Theorem
\ref{gqs}) that if $A$ is a ring such that the pro-vectorspace
\begin{equation}\label{barainfy}
H^{bar}_r((A\otimes\Q)^\infty)=H^{bar}_r(A^\infty)\otimes\Q=\{H^{bar}_r(A^n)\otimes\Q\}_n
\end{equation}
 is zero for all
$r\ge 0$, then $A$ is $\infty$-excisive for rational $K$-theory in
the sense of Section \ref{cqp}. As explained in the introduction,
for the particular case when $A=A^2$, this says the same as
Suslin-Wodzicki's theorem \eqref{tqs}. The strategy of proof
imitates that used by Suslin and Wodzicki in \cite{qs}. The three
main steps are Propositions \ref{siglex}, \ref{equiG} and
\ref{elbue}. Next we outline some of the contents of these propositions, and explain how the
main result of the section (\ref{gqs}) follows from them. In this summary, as well as in most of
the section, only homology with rational coefficients is considered, and the $\Q$ is dropped from
our notation; thus here $H_*(\ \ )=H_*(\ \ ,\Q)$. In \ref{siglex}, we consider the affine group
$\widetilde{GL}(A)=GL(A)\ltimes M_{1\infty}A$.
We show that if the map
\begin{equation}\label{piso}
H_r(GL(A^\infty))\to H_r( \widetilde{GL}(A^\infty))
\end{equation}
is an isomorphism for all $r$ and the same is true with $A^{op}$ substituted
for $A$, then $A$ is $\infty$-excisive for rational $K$-theory. In \ref{equiG} it is shown
that the condition that \eqref{piso} be an isomorphism is equivalent to two other conditions.
One of these (\ref{equiG} (b)) involves the space $\cup_{\tau,n}BT_n^\tau\subset BGL$, union over all $n$
and all finitely supported partial orders $\tau$ of the classifying spaces of the triangular
groups $T_n^\tau\subset GL_n$ defined by the orders.
The condition is that
\begin{equation}\label{qiso}
 \{H_r(\bigcup_{n,\tau}BT_n^\tau(A^m))\}_m\to
\{H_r(\bigcup_{n,\tau}B{{\widetilde{T}_n}^\tau}(A^m))\}_m
\end{equation}
be an isomorphism for all $r$. In particular \eqref{piso} is an isomorphism if \eqref{qiso} is.
In \ref{elbue}(b) we show that the condition that
\begin{equation}\label{riso}
H^{bar}_r(A^\infty)\otimes\Q=0
\end{equation}
for all $r$ implies \eqref{qiso} is an isomorphism. Thus \eqref{riso} implies \eqref{piso} is an
isomorphism. But we prove
in \ref{aaop} that $H^{bar}$ does not change if we replace an algebra by its opposite.
Hence \eqref{riso} implies that \eqref{piso} is an isomorphism not only for $A$ but also for $A^{op}$,
whence $\eqref{riso}$ implies excision in rational $K$-theory, by \ref{siglex}.$\qed$
\bigskip

Recall that if $A\in\rings$, its general linear group is defined
as
$$
GL(A):=\ker(GL(\tilde{A})\to GL(\Z)).
$$
Here the map $GL(\tilde{A})\to GL(\Z)$ is that induced by
$\tilde{A}\to \Z$, $(a,n)\mapsto n$. We remark that if $A$ happens to be unital, then $GL(A)$
as defined above is naturally isomorphic to the usual $GL(A)$. More generally, if $G:\as1\to (({\rm Groups}))$
is any functor which preserves finite products and $A$ is unital then
$$
\ker(G(\tilde{A})\to G(\Z))=\ker(G(A)\times G(\Z)\to G(\Z))=G(A).
$$
We consider the affine group
\begin{equation}\label{affine}
\widetilde{GL}(A)=GL(A)\ltimes M_{1\infty}A.
\end{equation}

\begin{lem}\label{redu}
Let $A\in \rings$, and assume that for every diagram \eqref{diag}
the natural map of pro-vectorspaces
$$
\{K^\Q_n(B:A^m)\}_m\lra \{K^\Q_n(C:A^m)\}_{m}
$$
is an isomorphism for $n\ge 1$. Then $A$ is $\infty$-excisive for
$\K$.
\end{lem}
\begin{proof} There is a map of exact sequences
$$
\xymatrix{
{\lim_m}^1K^\Q_{n+1}(B:A^m)\ar[d]\ar@{^{(}->}[r]
&\pi_n{\hat{\mathbf K}}^\Q(B:A^\infty)\ar[d]\ar@{->>}[r]&\lim_mK^\Q_n(B:A^m)\ar[d]\\
{\lim_m}^1K^\Q_{n+1}(C:A^m)\ar@{^{(}->}[r]
&\pi_n{\hat{\mathbf K}}^\Q(C:A^\infty)\ar@{->>}[r]&\lim_mK^\Q_n(C:A^m)}
$$
The hypothesis implies that the vertical arrows on both extremes
are isomorphisms for $n\ge 1$ and that the leftmost vertical arrow
is an isomorphism for $n=0$. Because $K$-theory is excisive in
dimension $\le 0$, the map
$$
K^\Q_n(B:A^m)\lra K^\Q_n(C:A^m)
$$
is an isomorphism for all $m$ and all $n\le 0$. We conclude that, under the
hypothesis of the lemma, the vertical map in the middle of
the diagram above is an isomorphism for all $n\in\Z$.
\end{proof}

\begin{lem}\label{FBA} Let $R$ be a unital ring, $A\triqui R$ an ideal and $\overline{GL}(R/A)$ the
image of $GL(R)$ in $GL(R/A)$. Write $F(R:A)$ for the fiber of the fibration $K^\Q(R)\to
\Q_{\infty}{B\overline{GL}}(R/A)$, and $F'(R:A):=\hofi(K^\Q(A)\to
K^\Q(R))$. Then the canonical map $F(R:A)\to F'(R:A)$ induces a weak equivalence onto the connected
component of $F'(R:A)$.
\end{lem}
\begin{proof}
We have a map of homotopy fibration sequences
\begin{equation}\label{fba}
\xymatrix{F(R:A)\ar[r]\ar[d]&K^\Q(R)\ar[r]\ar[d]_1&\Q_\infty(B\overline{GL}(R/A))\ar[d]\\
          F'(R:A)\ar[r]&K^\Q(R)\ar[r]&K^\Q(R/A)}
\end{equation}
From the long exact sequence of homotopy groups of the top fibration, we get that $F(R:A)$
is connected. Hence it suffices to prove
$\pi_nF(R:A)\to \pi_nF'(R:A)$ is an isomorphism for $n\ge 1$. In turn, by the five lemma applied to the map
of long exact sequences of homotopy groups associated to \eqref{fba}, we are further reduced to
proving that the right vertical arrow induces a isomorphism at the level of $\pi_n$ for $n\ge 2$ and
an injection in $\pi_1$. We have the identities
\begin{multline*}
[GL(R/A),GL(R/A)]=E(R/A)=\rm{Image}(E(R)\to \overline{GL}(R/A))\\
=\rm{Image}([GL(R),GL(R)]\to GL(R/A))=[\overline{GL}(R/A),\overline{GL}(R/A)]
\end{multline*}
Put $\overline{K}_1(R/A):=\rm{Image}(K_1(R)\to K_1(R/A))$. We have a map of fibrations
\[
\xymatrix{BE(R/A)\ar[d]_1\ar[r]&B\overline{GL}(R/A)\ar[d]_{\iota}\ar[r]&B\overline{K}_1(R/A)\ar[d]\\
BE(R/A)\ar[r]&BGL(R/A)\ar[r]& BK_1(R/A)}
\]
Apply $\Q_\infty$ to the diagram above. Because $K_1(R/A)$ is commutative,
\[\Q_\infty B\overline{K}_1(R/A)\iso B(\Q\otimes \overline{K}_1(R/A))\text{ and }\Q_\infty B K_1(R/A)\iso B(\Q\otimes K_1(R/A)).\]
Thus $\pi_n\Q_\infty(\iota)$ is an isomorphism for $n\ge 2$ and the inclusion $\Q\otimes \overline{K}_1(R/A)\subset \Q\otimes K_1(R/A)$
for $n=1$.
\end{proof}

\begin{notation} From now until the end of the current section, we shall mostly consider
homology with rational
coefficients, hence we shall drop the coefficients from our
notation, except when both rational and integral homology appear,
in which case the coefficients will be emphasized. Thus if $X$ is
a space, $H_*(X)$ will mean $H_*(X,\Q)$, and if $G$ is a group,
$$
H_*(G)=H_*(BG)=H_*(G,\Q).
$$
\end{notation}

\begin{prop}\label{siglex}
Let $A\in\rings$ and $\widetilde{GL}(A)$ be as in \eqref{affine}.
Assume that the inclusion $GL(A)\subset \widetilde{GL}(A)$ induces
an isomorphism of pro-abelian groups
\begin{equation}\label{ah}
H_r(GL(A^\infty))\iso H_r( \widetilde{GL}(A^\infty))
\end{equation}
for all $r\ge 1$, and that the same is true with $A^{op}$
substituted for $A$. Then $A$ is $\infty$-excisive for $\K$.
\end{prop}
\begin{proof}
Let $\iota:GL(A)\subset\widetilde{GL}(A)$ be the inclusion and
$p:\widetilde{GL}(A)\to GL(A)$ the projection. Because $p\iota=1$,
$\iota p$ induces the identity in pro-homology. In terms of
transition maps, this means
that for every $n$ and every $r$ there exists $k=k(n,r)\ge n$ such that for
the
$\iota p\sigma_{k,n}=\sigma_{k,n}$ on $H_r(\widetilde{GL}(A^k))$ (notations
are as in \ref{nota}).
Since $A^{op}$ also satisfies this by hypothesis, and since both $p$ and
$\iota$ are level maps, and moreover, by naturality, the $\sigma$ commute
with the action of $Gl(\Z)$, the argument of \cite[Prop. 1.5]{qs},
shows that
\begin{equation}\label{zfix}
\sigma_{k,n}(H_r(GL(A^k))) \subset (H_r(GL(A^n)))^{GL(\Z)}.
\end{equation}
Next the argument of \cite[Cor. 1.6]{qs} proves that if
$A\triangleleft B$ is an embedding as an ideal of a unital ring
$B$ then \eqref{zfix} holds with $GL(B)$ subsituted for $GL (\Z)$
\begin{equation}\label{rfix}
\sigma_{k,n}(H_r(GL(A^k))) \subset (H_r(GL(A^n)))^{GL(B)}.
\end{equation}
Now proceed as in the proof of \cite[1.7]{qs}, using the
levelwise Serre spectral sequence and then \ref{trivact} and
\ref{compa}, to conclude that the map $BGL(A^\infty)\to F(B:A^\infty)$ induces
an isomorphism in rational pro-homology.  In particular if $B\to C$ is as in
\eqref{diag}, then the map $\alpha:F(B:A^\infty)\to F(C:A^\infty)$ is a rational pro-homology equivalence, by the argument of the proof of \cite[1.7]{qs}.
Because by \ref{FBA} both $F(B:A^\infty)$ and $F(C:A^\infty)$ are towers of connected spaces each of which is a homotopy loopspace whose homotopy
groups are $\Q$-vectorspaces, it follows from Prop. \ref{minomu} that $\alpha$ is a pro-homotopy equivalence. By the previous lemmas, this shows that $A$ is $\infty$-excisive.
\end{proof}

\begin{notation}\label{notatriangu}
In the next proposition, $E(A)\subset GL(A)$ is the elementary
subgroup, $T_n^\tau(A)\subset GL_n(A)$ is the triangular subgroup
defined by the partial order $\tau$ of the set $\{1,\dots,n\}$,
and $\cup_{\tau,n}BT_n^\tau\subset BGL$ is the union over all $n$
and all finitely supported partial orders $\tau$. The definition
of each of the groups just mentioned for nonunital $A$ can be
found in \cite{qs}. If $G\subset GL(A)$ is a subgroup, we write
$$
\widetilde{G}:=G \ltimes M_{1\infty}A.
$$
\end{notation}

\begin{prop}\label{equiG}
Let $A\in\rings$, $\iota:GL(A)\to\widetilde{GL}(A)$ the inclusion.
The following are equivalent.
\item{(a)} The map $\iota$ induces a pro-isomorphism
$$
H_r(GL(A^\infty))\to H_r(\widetilde{GL}(A^\infty))\qquad (\forall
r\ge 0).
$$
\item{(b)} The map $\iota$ induces a pro-isomorphism
$$
H_r(E(A^\infty))\to H_r(\widetilde{E}(A^\infty))\qquad (\forall
r\ge 0).
$$
\item{(c)} The map $\iota$ induces a pro-isomorphism
\begin{equation}\label{equiT}
 \{H_r(\bigcup_{n,\tau}BT_n^\tau(A^m))\}_m\to
\{H_r(\bigcup_{n,\tau}B{{\widetilde{T}_n}^\tau}(A^m))\}_m\qquad
(\forall r\ge 0).
\end{equation}
\end{prop}
\begin{proof}
Let $C(A)=[GL(A),GL(A)]$. We have a natural map of short exact
sequences
$$
\begin{CD}
C(A) @>>> GL(A) @>>> GL_{ab}(A)\\
@VVV @VVV @VV\cong V\\
\tilde{C}(A)@>>>\widetilde{GL}(A) @>>>\widetilde{GL}(A)/\tilde{C}(A)\\
\end{CD}
$$
This map induces a level map of pro-Lyndon-Serre spectral sequences
$$
\xymatrix{E^2_{p,q}=\{H_p(GL_{ab}(A^n),H_q(C(A^n)))\}_n\ar@{=>}[r]\ar[d]&
H_{p+q}(GL(A^n))\ar[d]\\
E^2_{p,q}=\{H_p(GL_{ab}(A^n),H_q(\tilde{C}(A^n)))\}_n\ar@{=>}[r]&
H_{p+q}(\widetilde{GL}(A^n))}
$$
By the proof of \cite[Cor. 1.14]{qs},
(see formula (12)),
$E(A^{2^n})$ is closed under conjugation by elements of
$GL(A^{2^{n+1}})$. The inclusions
$$
C(A^{2^{n+1}})\subset E(A^{2^n}),\quad E(A^{2^{n+1}})\subset
[E(A^{2^n}),E(A^{2^n})]\subset C(A^{2^n})
$$
induce, for each $n\ge 2$, two homomorphisms of pairs (see the
Appendix for a definition)
$$
(GL(A^{2^{n+1}}),C(A^{2^{n+1}}))\to
(GL(A^{2^{n+1}}),E(A^{2^n}))\to (GL(A^{2^{n-1}}),C(A^{2^{n-1}})).
$$
These give mutually inverse pro-isomorphisms
$$
\{(GL(A^n),C(A^n))\}_n\overset{\sim}\longleftrightarrow\{(GL(A^{2^{n+1}}),E(A^{2^n}))\}_n.
$$
Taking homology we obtain a pro-isomorphism
$$
\{H_p(GL_{ab}(A^n),H_q(C(A^n)))\}_n\cong
\{H_p(GL_{ab}(A^{2^{n+1}}),H_q(E(A^{2^n})))\}_n.
$$
 Similarly,
$$
\{H_p(GL_{ab}(A^{n}),H_q(\tilde{C}(A^{n})))\}_n\cong
\{H_p(GL_{ab}(A^{2^{n+1}}), H_q(\tilde{E}(A^{2^n})))\}_n.
$$
The proof that (a)$\iff$(b) now follows as in the proof of \cite[2.10]{qs},
using \cite[1.14]{qs} and
\ref{trivact}, \ref{elexa} and \ref{compa}. Similarly, that (b)$\iff$(c) follows from
\cite[2.9 and 2.7]{qs}, and
\ref{trivact} and \ref{compa}.
\end{proof}

\begin{notation}
In the following lemma and below, $\gl_n(A)$ is the general linear
Lie algebra, $\pri_n(A)=[\gl_n(A),\gl_n(A)]$,
$\rip_n(A)=\ker($Trace$:\gl_n(A)\to A/[A,A])$, and  $\e_n(A)\subset\gl_n(A)$ is the Lie subalgebra generated by
the matrices $e_{i,j}(a)$ ($i\ne j$, $a\in A$).
\end{notation}

\begin{lem}\label{identi}
Let $A\in\ass(\Q)$, $n\ge 3$. Then
\begin{align*}\nonumber
\e_n(A)=\bigoplus_{1\le i\ne j\le n}e_{i,j}(A)\oplus\bigoplus_{i=1}^{n-1}
(e_{i,i}-e_{{i+1},{i+1}})(A^2)\oplus e_{1,1}([A,A])\nonumber\\
\rip_n(A)=\bigoplus_{1\le i\ne j\le n}e_{i,j}(A)\oplus\bigoplus_{i=1}^{n-1}
(e_{i,i}-e_{{i+1},{i+1}})(A)\oplus e_{1,1}([A,A])\nonumber\\
\pri_n(A)=\bigoplus_{1\le i\ne j\le n}e_{i,j}(A^2)\oplus
\bigoplus_{i=1}^{n-1}(e_{i,i}-e_{i+1,i+1})(A^2)\oplus
e_{1,1}([A,A]).\nonumber\\
\end{align*}

\end{lem}
\begin{proof} Straightforward.\end{proof}
\begin{cor}\label{fil}
\begin{align*}
\rip_n(A^2)\subset  \pri_n(A)&\subset \rip_n(A)\nonumber\\
\rip_n(A^2)\subset  \e_n(A)&\subset \rip_n(A).\qed\nonumber
\end{align*}
\end{cor}

\smallskip

\begin{notation}
In the following proposition, $\tri_n^\tau(A)\subset\gl_n(A)$ is
the triangular Lie subalgebra defined by the partial order $\tau$
of $\{1,\dots,n\}$ (see \cite{qs}). We put
$$
\gl_{n+1}(A)\supset\widetilde{\gl}_n(A)=\{\begin{pmatrix}g&v\\ 0&0\end{pmatrix}
:g\in\gl_n(A), v\in M_{n,1}(A)\}
$$
and for $\g=\rip,\pri,\e,\dots$
$$
\tilde{\g}_n(A)=\g_{n+1}(A)\cap \widetilde{\gl}_n(A).
$$
If $\g_n(A)$ is either $\gl_n(A)$ or one of the Lie algebras
above, then $\g(A)$ with no subscript stands for $\cup_n\g_n(A)$.
In the proposition below, $A$ is a $\Q$-algebra and  Lie algebra
homology is taken with coefficients in $\Q$. We write
$$
H_*(\g):=H_*(\g,\Q)
$$
for all Lie algebras $\g$ appearing below. The groups $H_*(\g)$
are the homology of the the Chevalley-Eilenberg complex $C_*(\g)$ for Lie algebra homology
over $\Q$.
(\cite[10.1.3]{lod}). This complex appears in item $(c)$ of the
proposition below.
\end{notation}
\begin{prop}\label{equig}
Let $A\in\ass(\Q)$, $\iota:\gl(A)\to\widetilde{\gl}(A)$ the
inclusion. The following are equivalent.
\item{(a)} The map $\iota$ induces a pro-isomorphism
$$
H_r(\gl(A^\infty))\to H_r(\widetilde{\gl}(A^\infty))\qquad (r\ge
0).
$$
\item{(b)} The map $\iota$ induces a pro-isomorphism
$$
H_r(\e(A^\infty))\to H_r(\widetilde{\e}(A^\infty))\qquad (r\ge 0).
$$
\item{(c)} The map $\iota$ induces a pro-isomorphism
$$
\{H_r(\sum_{n,\tau}C_*(\tri_n^\tau(A^m))\}_m\to
\{H_r(\sum_{n,\tau}C_*({\widetilde{\tri}_n}^\tau(A^m))\}_m\qquad
(r\ge 0).
$$
\end{prop}
\begin{proof}
As in \cite[formula (34)]{qs},
we shall write $\nu(A)$ and
$\tilde{\nu}(A)$ for the Volodin chain complexes. To start, mimic
the proof of \cite[2.9]{qs}, to show that the actions of $\e(A^2)$
on $H_*(\nu(A))$ and of $\tilde{\e}(A^2)$ on $H_*(\tilde{\nu}(A))$
are trivial. Then use \ref{identi} and imitate the proof of
\cite[1.14]{qs}, to get that the actions of $\gl(A^2)$ on
$H_*(\rip(A))$, of $\gl(A^2)$ on $H_*(\e(A))$, of
$\widetilde{\gl}(A^2)$ on $H_*(\widetilde{\rip}(A))$ and of
$\widetilde{\gl}(A^2)$ on $H_*(\tilde{\e}(A))$ are trivial. Then
follow as in the proof of \ref{equiG}, substituting the
Hochschild-Serre spectral sequence for that of Lyndon-Serre, and
\cite[4.9]{qs} for \cite[2.7]{qs}.
\end{proof}

\begin{notation}
If $R$ is a ring, we write $R_1$ for the ring of $2\times 2$ matrices
$$
M_2R\supset R_1:=\{\begin{pmatrix}a &b\\ 0& 0\end{pmatrix}: a,b\in
R\}.
$$
A straightforward calculation shows that the construction $R\mapsto R_1$ commutes with powers
$$
(R^n)_1=(R_1)^n.
$$
In particular if $A$ is any ring, then the pro-ring $A_1^\infty=\{A_1^n\}_n$ is unambiguously defined.

In the proof of the lemma below, the {\it bar complex} $C_*^{bar}(A)$
of a ring $A$ is considered. Recall from \cite{wod}
that $C^{bar}_n(A)=0$ if $n\le 0$ and
$C^{bar}_n=A^{\otimes n}$ if $n\ge 1$, with boundary operator
given by
$$
b'(a_1\otimes \dots\otimes a_n)=\sum_{i=1}^{n-1}(-1)^{i+1}a_1\otimes\dots \otimes
a_i a_{i+1}\otimes\dots
\otimes a_n.
$$
By definition, the bar homology of $A$ is the homology of $C^{bar}_*(A)$.

\end{notation}

\begin{lem}\label{100pages}
Let $A\in\ass(\Q)$. Assume $H^{bar}_r(A^\infty)=0$ ($r\ge 0$). Then
\item{(a)}
\begin{equation}\label{bar1nul}
H^{bar}_r(A_1^\infty)=0 \qquad (r\ge 0)
\end{equation}
\item{(b)} The inclusion $A\subset A_1$ induces a pro-isomorphism
\begin{equation}\label{hca1}
HC_r(A^\infty)\cong HC_r(A_1^\infty)\qquad (r\ge 0)
\end{equation}
\end{lem}
\begin{proof}
(a) This is straightforward from the proof of a particular case of
\cite[9.1]{wod}. We recall
the argument for the reader's convenience. If $R$ is any $\Q$- algebra, we write $N(R)$ for the following ideal of $R_1$
$$
N(R):=\{\begin{pmatrix}0 &b\\ 0& 0\end{pmatrix}: b\in R\}.
$$
We may regard $N(R)$ as an $R$-bimodule; note that $R_1$ is isomorphic to the semidirect
product $R\ltimes N(R)$. Hence the proof of \cite[9.1]{wod} applies to $R_1$, and shows in particular
(see the first displayed formula on page 622 of {\it loc.cit.}) that the bar complex
has a direct sum decomposition
$$
C^{bar}(R_1)=C^{bar}(R)\oplus \bigoplus_{l=1}^\infty C^{bar}(R_1,l)
$$
By \cite[9.2]{wod}, for each $l\ge 1$ there is a first quadrant spectral sequence converging to the homology
of $C^{bar}(R_1,l)$ of which the zero term is of the form
\begin{equation}\label{speceasy}
(E_{p*}^0(R_1,l),d^0)=C^{bar}(R)\otimes V(R,l)[l-1]
\end{equation}
where $V(R,l)$ is some vectorspace depending on $R$ and $l$. Applying all this levelwise to
$A_1^\infty$, and noting that the hypothesis that the bar pro-homology of $A^\infty$ vanishes implies
\begin{equation}\label{barvnul}
\{H^{bar}_r(A^n)\otimes V_n\}=0\qquad (r\ge 0)
\end{equation}
for every pro-vectorspace $V$, proves (a).

\smallskip

(b) It suffices to show that the inclusion $A\subset A_1$ induces a pro-isomorphism
$HH_r(A^\infty)\cong HH_r(A_1^\infty)$. To prove this we use the argument of the proof
of \cite[11.1]{wod}. It is shown in {\it loc.cit.} that for $R\in \ass(\Q)$ the Hochschild
complex decomposes as $C(R_1)=C(R)\oplus\bigoplus_{l=1}^\infty C(R,l)$ (see formula
(50) in {\it loc.cit.}) and that again for each $l\ge 1$ there is a first quadrant spectral sequence converging to $H_*(C(R,l))$
of the form \eqref{speceasy} (see \cite[11.2]{wod}). As in (a), applying this levelwise to
$A_1^\infty$, we obtain (b).
\end{proof}

\begin{prop}\label{elbue}
\item{(a)} Let $A\in\ass(\Q)$. Assume
\begin{equation}\label{barnul}
H^{bar}_r(A^\infty)=0 \qquad (r\ge 0).
\end{equation}
Then the natural inclusions $T_n(A)\subset\tilde{T}_n(A)$ ($n\ge
1$) induce a pro-isomorphism
$$
\{H_r(\bigcup_{n,\tau}BT_n^\tau(A^m),\Z)\}_m\liso
\{H_r(\bigcup_{n,\tau}B\tilde{T}_n^\tau(A^m),\Z)\}_m\qquad (r\ge 0).
$$
\item{(b)} Let $A\in\rings$. Assume
\begin{equation*}
H^{bar}_r(A^\infty)\otimes\Q=0\qquad (r\ge 0).
\end{equation*}
Then the natural inclusions $T_n(A)\subset\tilde{T}_n(A)$ ($n\ge
1$) induce a pro-isomorphism
$$
\{H_r(\bigcup_{n,\tau}BT_n^\tau(A^m),\Q)\}_m\liso
\{H_r(\bigcup_{n,\tau}B\tilde{T}_n^\tau(A^m),\Q)\}_m\qquad (r\ge
0).
$$
\end{prop}
\begin{proof}
(a)
According to \cite[Cor. 5.14]{qs}, for any
$m>0$ we have canonical isomorphisms
\begin{equation}\label{pegate}
\begin{split}
\tilde{H}_*(\bigcup_{n,\tau}BT^\tau_n(A^m),\Z)\iso
\tilde{H}_*(\sum_{n,\tau}C_*(t^\tau_n(A^m)))\\
\tilde{H}_*(\bigcup_{n,\tau}B\tilde{T}^\tau_n(A^m),\Z)\iso
\tilde{H}_*(\sum_{n,\tau}C_*(\tilde{t}^\tau_n(A^m))).\\
\end{split}
\end{equation}
Thus to establish the required pro-isomorphism it suffices to show
that the canonical map
$$
\{H_*(\sum_{n,\tau}C_*(t^\tau_n(A^m)))\}_m\to
\{H_*(\sum_{n,\tau}C_*(\tilde{t}_n^\tau(A^m)))\}_m
$$
is a pro-isomorphism, which according to Proposition \ref{equig}
is equivalent to verifying that the canonical map
$H_*(\gl(A^\infty))\to H_*(\widetilde{\gl}(A^\infty))$ is a
pro-isomorphism. Since the pro-Lie algebra
$\widetilde{\gl}(A^\infty)$ is a retract in $\gl(A_1^\infty)$, it is
enough to show that the embedding
\begin{equation}\label{gluglu}
\gl(A^\infty)\to \gl(A_1^\infty)
\end{equation}
induces an isomorphism in pro-homology. By \eqref{barnul} and \ref{100pages},
\eqref{bar1nul} holds for $A$. But according to \cite[Thm. 4.2]{hanwi} if
$P$ is any pro-algebra such that $H^{bar}_l(P)=0$ for all $l$, then
\[
H_r(\gl(P))\cong (\wedge (HC(P)[-1]))_r
\]
for all $r$. Here $HC(P)=\{\bigoplus_{n\ge 0}HC_n(P_m)\}_m$; the $[-1]$ indicates a levelwise degree shift. Thus
\begin{equation}\label{hglhc}
\begin{split}
H_r(\gl(A^\infty))\cong (\wedge (HC(A^\infty)[-1]))_r\\
H_r(\gl(A_1^\infty))\cong (\wedge ( HC(A_1^\infty)[-1]))_r.\\
\end{split}
\end{equation}
By virtue of \eqref{hca1}, this implies that \eqref{gluglu} induces an
isomorphism in pro-homology, as wanted.
Part (b) follows from part (a) and \cite[Cor. 5.19]{qs}.
\end{proof}

\begin{lem}\label{aaop}
The map
$$
y:C^{bar}_*(A)\to C^{bar}_*(A^{op}),\quad y(a_1\otimes\dots\otimes
a_n)=(-1)^{\frac{n(n-1)}{2}}a_n\otimes\dots\otimes a_1
$$
is an isomorphism of chain complexes.
\end{lem}
\begin{proof}
Straightforward.
\end{proof}

\begin{thm}\label{gqs}
Let $A$ be a ring. Assume that $H_r^{bar}(A^\infty)\otimes\Q=0$
($r\ge 0$). Then $A$ is $\infty$-excisive (in the sense of Section
\ref{cqp}) for rational $K$-theory.
\end{thm}
\begin{proof}
By \ref{siglex} and \ref{aaop}, to prove $A$ is $\infty$-excisive
for $\K$ it suffices to show that for all $r\ge 0$, the map
\eqref{ah} is an isomorphism. By \ref{equiG} this is equivalent to
showing that \eqref{equiT} is an isomorphism, which is in fact the
case, as proved by Proposition \ref{elbue}.
\end{proof}

\section{Proof of the main theorem}\label{qsp}

\noindent{\sc Summary.} Here we carry out part (d) of the sketch
of the introduction and complete the proof of the Main theorem
(4.2). We prove first (Lemma \ref{chum}) that for $\Q$-algebras
$A$ and $B$ the $K$-theory obstruction groups $K_*(A,B:I)$ are
$\Q$-vectorspaces. As explained in the introduction, after this
lemma and Theorems \ref{fund}, \ref{agree} and \ref{gqs}, to
finish the proof it remains only to show that if $A$ is an ideal
of a unital free algebra then $H^{bar}_r(A^\infty)=0$ ($r\ge 0$)
and $A$ is $\infty$-excisive for $HC$. We explain how to derive
this from what is proved in \cite{val}, \cite{cq} and
\cite{soval}.$\qed$

\bigskip

\begin{lem}\label{chum} Let $A\to B$
be a homomorphism of $\Q$-algebras carrying an ideal
$I\triangleleft A$ isomorphically onto an ideal of $B$. Then the
map
$$
K_*(A,B:I)\lra K_*^\Q(A,B:I)
$$
is an isomorphism.
\end{lem}
\begin{proof} By \ref{exnuni1}, we may assume that $A\to B$ is a homomorphism of unital algebras.
It suffices to show that the group $K_*(A,B:I)$ is uniquely
divisible. Let $m\in \Z$; there is a long exact sequence
\begin{multline}
K_{n+1}(A,B:I,\Z/m)\to K_n(A,B:I) \overset{\cdot m}\to\\
K_n(A,B:I) \to K_n(A,B:I,\Z/m)
\end{multline}
By \cite{chu} the groups in both extremes are zero.
\end{proof}

\noindent{\bf 4.2. Proof of the Main Theorem.} By the previous
lemma, it suffices to consider rational $K$-theory. By Corollary
\ref{exnuni1}, Theorems \ref{fund} and \ref{agree}, and Remarks
\ref{exnuni3} and \ref{chetaununi}, this will
follow if we show that if $A$ is an ideal of a unital free ring,
then $A$ is $\infty$-excisive for both $K^\Q$ and $HN^\Q$. But if
$H^{bar}_r(A^\infty)\otimes\Q=0$ for $r\ge 0$, then $A$ is
$\infty$-excisive for $K^\Q$ by Theorem \ref{gqs}. As explained in
the sketch of the proof of the Main Theorem in the introduction
(part (d)), it suffices to prove that if $A$ is an ideal of a free
unital $\Q$-algebra then
\begin{equation}\label{barzero}
H^{bar}_r(A^\infty)=0\qquad (r\ge 0)
\end{equation}
and that $A$ is $\infty$-excisive for $HC$. Both the latter
property and \eqref{barzero} are well-known, and, when properly
restated in terms of cofibrations and weak equivalences, hold with
much greater generality than needed here; for example in
\cite[Thms. 5.2 and 6.3 ]{val},
they are proved for pro-algebra
objects in a linear category with a tensor product. We shall
explain presently how a more direct proof of \eqref{barzero} in
our specific situation can be extracted from \cite{cq} and
\cite{soval}. First of all we observe that even if \cite{soval} is
written for the topological algebra setting the algebraic case
works in exactly the same way. From \cite[Prop. 4.2]{cq},
$A$ is
quasi-unital in the sense of \cite[Def. 4.6]{soval}. Proposition
4.7 of \cite{soval} states that, for quasi-unital $A$, the bar
pro-complex $C^{bar}_*(A^\infty)$ is weakly contractible. Actually
the proof of \cite[4.7]{soval} shows that for each $n$ the map
\begin{equation}\label{transfer}
C^{bar}_*(A^{2^{n+1}})\to C^{bar}_*(A^{2^n})
\end{equation} is
null-homotopic; \eqref{barzero} is immediate from this. Next, to
prove $A$ is $\infty$-excisive for $HC$, it suffices, by virtue of
the version of the $SBI$-sequence involving $HC$ and Hochschild
homology $HH$ (\cite[2.5.8]{lod}), to show that $A$ is
$\infty$-excisive for $HH$. This means that if $B$ is a unital
algebra and $A\triqui B$, then for the Hochschild complex
$\Omega^*$ and the Dold-Kan functor $S$, the map
$$
\holi_nS\Omega^*(\tilde{A}:A^n)\to \holi_nS\Omega^*(B:A^n)
$$
is an equivalence. In view of the exact sequence
$$
0\to {\lim_n}^1\pi_{r+1}X^n\to \pi_r\holi_nX^n\to\lim_n\pi_r
X^n\to 0
$$
it suffices to show that for each $r$ the map
\begin{equation}\label{flechi}
HH_r(\tilde{A}:A^\infty)\to HH_r(B:A^\infty)
\end{equation}
is an isomorphism of pro-vectorspaces. One can deduce this
directly from \eqref{barzero} by mimicking a proof of Wodzicki's
theorem, such as for example the brief proof given in \cite{gg}.
Alternatively one can deduce it from the results of \cite{soval}
and \cite{val} as follows.  We have proved above that
$C^{bar}(A^\infty)$ is weakly contractible, which by \cite[Thm. 5.2]{val},
(or \cite[4.4]{soval}) implies that
$\Omega^*(\tilde{A}:A^\infty)\to\Omega^*(B:A^\infty)$ is a weak
equivalence (a finite weak equivalence in the notation of
\cite{soval}) of pro-complexes. Hence \eqref{flechi} is an
isomorphism by \cite[Lem. 2.1.1]{val}.
\qed

\begin{rem}
If $A\in\ass(\Q)$ then by the Main Theorem and Remark \ref{chetauq},
$$A\mapsto\Hy(inf(\ass(\Q)\downarrow A),\mathbf{K})$$ is excisive.
\end{rem}

\appendix

\section{Pro-spaces}

\noindent{\sc Summary.} We introduce the notation for pro-objects
used throughout the paper, and prove some technical results on pro-spaces
which are used in Section \ref{sgqs}.
Two of these concern the Serre spectral sequence
(\ref{trivact} and \ref{compa}), and another rational pro-homotopy
and pro-homology \eqref{minomu}.$\qed$

\bigskip

\begin{notation}\label{nota}
If $\C$ is a category, we write pro-$\C$ for the category of all
inverse systems
$$
\{\sigma_{n+1}:C_{n+1}\to C_n:n\in\Z_{\ge 1}\}
$$
of objects of $\C$. The set of homomorphisms of two pro-objects
$C$ and $D$ in pro$-\C$ is by definition
$$
\hom_{\pro-\C}(C,D)=\lim_n\coli_m\hom_\C(C_m,D_n)
$$
Thus a map $C\to D\in\pro-\C$ is an equivalence class of maps of
inverse systems
$$
f:\{C_{k(n)}\}_n\to \{D_n\}_n
$$
where $k:\Z_{\ge 1}\to \Z_{\ge 1}$ satisfies $k(n)\ge n$ for all
$n$. In the particular case when $k$ is the identity we say that
$f$ (and also its equivalence class) is a {\it level map}. We use
the letter $\sigma$ (often with no subscript) to mean the
transition maps of any pro-objects appearing in this paper. If
$k\ge n$, we write $\sigma_{k,n}$ for the composite
$$\sigma_{k,n}:=\sigma_{n+1}\circ\dots\circ\sigma_k.$$
Note that our $\pro-\C$ is equivalent to the category of countably
indexed pro-objects considered by Cuntz and Quillen (\cite{cq});
both are much smaller than that of Artin-Mazur (\cite{am}).

\smallskip

If $(G,L)$ is a pair consisting of a group $G$ and a left
$G$-module $L$ then by a morphism from $(G,L)$ to another such
pair $(H,M)$ we mean a group homomorphism $f:G\to H$ together with
a homomorphism of abelian groups $L\to M$ --which we shall also
call $f$-- such that
$$
f(gl)=f(g)f(l)\qquad(\forall g\in G, l\in L).
$$
We write $C_*(G,M)$ for the bar complex of $G$ with coefficients
in $M$.
\end{notation}

\begin{lem}\label{trivact}

Let $(G,M)=\{\sigma_{n+1}:(G_{n+1},M_{n+1})\to (G_{n},M_n):n\ge
1\}$ be a sequence of pairs consisting at each level $n$ of a
group $G_n$ and a left $G_n$-module $M_n$, together with
transition maps $\sigma_n$ which are morphisms of such pairs in
the sense of \ref{nota}. Assume that
$$
(\forall n)(\exists k=k(n)>n)\text{ such that }(\forall g\in G_k,
x\in M_k) \ \ \sigma_{k,n}(gx)=\sigma_{k,n}(x).
$$
Then there is a natural short exact sequence pro-graded abelian
groups
\begin{equation}\label{ses}
0\to \{ H_*(G_n,\Z)\otimes_{\Z}M_n\}\to
\{H_*(G_n,M_n)\}\to\{\tor(H_{*-1}(G_n,\Z), M_n)\}\to 0.
\end{equation}
\end{lem}
\begin{proof}
In the particular case when each $G_n$ acts trivially on
each $M_n$, we have an inverse system of levelwise exact sequences
\begin{equation*}
0\to H_*(G_{n},\Z)\otimes M_n\to H_*(G_{n},M_n) \to \tor
(H_{*-1}(G_{n},\Z),M_n)\to 0.
\end{equation*}
This gives \eqref{ses} in this particular case. Next we observe
that applying the functor $(G,M)\to C_*(G,M)$ maps isomorphic
pro-pairs to isomorphic pro-complexes, which lead to isomorphic
pro-homology. Hence it suffices to show that any inverse system
$\{(G_n,M_n)\}$ satisfying the hypothesis of the lemma is
isomorphic to one as in the particular case considered above. The
inverse systems formed by the maps
$\sigma_{k^{n}(1)}:(G_{k^n(1)},M_{k^n(1)})\to (G_{n},M_{n})$ and
$1:(G_{k^{n}(1)},M_{k^n(1)})\to (G_{k^n(1)},M_{k^n(1)})$ represent
two mutually inverse isomorphisms
$$
\{(G_{n},M_{n})\}\cong \{(G_{k^n(1)},M_{k^n(1)})\}
$$
We may therefore assume that for all $n$, $k(n)=n+1$. For each $n$
consider $\sigma M_{n}$ as a $G_n$-module via $\sigma$, and let
$\iota:\sigma M_n\subset M_{n-1}$ be the inclusion. Then the
inverse systems of maps $(1,\sigma):(G_n,M_n)\to (G_n,\sigma M_n)$
and $(\sigma_n,\iota):(G_n,\sigma M_n)\to (G_{n-1},M_{n-1})$
represent mutually inverse isomorphisms
$$
\{(G_n,M_n)\}\cong\{(G_n,\sigma M_n)\}
$$
By hypothesis, the action of $G_n$ on $\sigma M_n$ is trivial.
\end{proof}

\begin{notation}\label{edge}
Let $\A$ be an abelian category, $(E^2_{p,q}, d_2)$ a first
quadrant spectral sequence in $\A$ of homological type, and
$F,G:\A\times\A\to\A$ two functors. We require that both functors
be {\it biadditive}, by which we mean that they be additive
separately on each variable. We say that the spectral sequence is
{\it $(F,G)$-dominated} if there are exact sequences
\begin{equation}
0\to F(E^2_{p,0},E^2_{0,q}) \to E^2_{p,q}\to
G(E^2_{p-1,0},E^2_{0,q})\to 0.
\end{equation}
If $E'^2$ is another $(F,G)$-dominated spectral sequence, then by
an $(F,G)$-homomor- phism $E^2\to E'^2$ we understand a family of
homomorphisms $f_{p,q}:E_{p,q}^2\to E'^2_{p,q}$ compatible with
differentials and such that for all $p$ and $q$, the obvious
diagram involving $f_{p,q}$, $F(f_{p,0},f_{0,q})$ and
$G(f_{p-1,0}, f_{0,q})$, commutes.
 \end{notation}

\begin{exa}\label{elexa}
Let $f:G\to G''$ be a homomorphism of inverse systems of groups
such that each $f_n:G_n\to G"_n$ is surjective; put $G'=\{\ker
f_n\}$. We have an inverse system of spectral sequences
\begin{equation}\label{prospec}
{}_nE^2_{p,q}=H_p(G''_n,H_q(G'_n,\Z))\Rightarrow H_{p+q}(G_n,\mathbb{Z})
\end{equation}
We regard \eqref{prospec} as a spectral sequence $E^*$ in the
category $\A=\pro-Ab$ of pro-abelian groups. Because each $_nE^*$
is located in the first quadrant, $E^*$ is convergent. Assume
further that the action of $G''$ on $H_*(G';\mathbb{Z})$ satisfies the
hypothesis of \ref{trivact}. Then by \ref{trivact} we have an
exact sequence
\begin{equation}
0\to \{{}_nE_{p,0}^2\otimes {}_nE_{0,q}^2\}\to
E^2_{p,q}\to\{\tor(_nE^2_{p-1,0},{}_nE^2_{0,q})\} \to 0
\end{equation}
Thus $E^2$ satisfies the requirements of \ref{edge} with
$$F(M,N)=\{M_n\otimes N_n\}\text{ and }
G(M,N)=\{\tor(M_n,N_n)\}.$$
\end{exa}

\begin{prop}\label{compa}
Let $E^2\to \E^2$ be a map of $(F,G)$-dominated first quadrant
spectral sequences in an abelian category $\A$. Let $H$ and $H'$
be filtered graded objects such that $E^*$ converges to $H$ and
${\E}^*$ to $H'$. Consider the maps
\begin{equation}\label{isos}
E^2_{*,0}\to \E^2_{*,0},\quad E^2_{0,*}\to \E^2_{0,*},\quad H_*\to
H'_*
\end{equation}
If any two of the maps in \eqref{isos} is an isomorphism, then so
is the third.
\end{prop}

\begin{proof} Zeeman's proof of his comparison theorem \cite{zee} proves
this statement.
\end{proof}

Recall that if $Z$ is a space and $k$ a field, then the diagonal map
$\Delta:Z\to Z\times Z$ induces a $k$-coalgebra structure on
$H(Z,k)=\bigoplus_nH_n(Z,k)$. Recall further that the {\it primitive part} of
$H(Z,k)$ is the graded vectorspace given in degree $n$ by
\begin{equation}\label{sofinito}
\primito_nH(Z,k)=\ker(\tilde{\Delta}:H_n(Z,k)\to \bigoplus_{0\le p\le n}H_p(Z,k)\otimes H_{n-p}(Z,k))
\end{equation}
where $\tilde{\Delta}(\xi)=\Delta(\xi)-\xi\otimes 1-1\otimes\xi$.
In the next lemma we consider the pro-extension of the functor
$\primito_nH(\ \ , k)$.

\begin{lem}\label{lemino}
Let $k$ be a field and $f:X\to Y$ a morphism of pro-spaces. Assume that
for each $n$ the map $H_n(f)\otimes k:H_n(X,k)\to H_n(Y,k)$ is an isomorphism
of pro-vectorspaces. Then the induced map
$\primito_nH(X,k)\to \primito_nH(Y,k)$
is an isomorphism for each $n$.
\end{lem}
\begin{proof} Immediate from the observation that for each fixed $n$, $\primito_nH(X,k)$
depends only on finitely many homology pro-vectorspaces $H_p(X,k)$, namely on those with
$p\le n$ (see \eqref{sofinito}).
\end{proof}
\begin{prop}\label{minomu}
Let $f:X\to Y$ be a pro-map of towers of connected fibrant pointed spaces.
Assume that for each
$m$ the spaces $X_m$, $Y_m$ are homotopy loopspaces and that
$H_n(f)\otimes \Q:H_n(X,\Q)\to H_n(Y,\Q)$
is an isomorphism of pro-vectorspaces for each $n$. Then
$\pi_n(f)\otimes\Q$ is a isomorphism for each $n$.
\end{prop}
\begin{proof} We have a commutative diagram
\[
\xymatrix{\pi_n(X)\otimes\Q\ar[d]\ar[r]&\pi_n(Y)\otimes\Q\ar[d]\\
          \primito_nH(X,\Q)\ar[r]&\primito_nH(Y,\Q)}
\]
where the horizontal maps are induced by $f$, and the vertical ones by the levelwise Hurewicz maps. In view
of the hypothesis on the $X_m$ and $Y_m$, the vertical maps are levelwise
isomorphisms. Because by hypothesis $H_p(f)\otimes\Q$ is an isomorphism for all $p$,
the bottom
row is an isomorphism as well, by Lemma \ref{lemino}. It follows
that the top row is an isomorphism.
\end{proof}


\begin{thebibliography}{10}
\bibitem{am} M. Artin, B. Mazur, {\it Etale homotopy}, Springer Lect. Notes Math. {\bf 100}
(1969).
\bibitem{bass} H. Bass, {\it Algebraic $K$-theory}. W. A.
Benjamin, New York Amsterdam 1968.
\bibitem{bk} A.K. Bousfield, D. Kan, {\it Homotopy limits, completions and
localizations}, Springer Lect. Notes Math. {\bf 304} (1972).
\bibitem{chino} G. Corti\~nas, {\it On the derived functor analogy
in the Cuntz-Quillen framework for cyclic homology} Algebra
Colloquim {\bf 5} (1998) 305-328.
\bibitem{crelle} G. Corti\~nas, {\it Infinitesimal $K$-theory}, J. reine
angew. Math. {\bf 503} (1998) 129-160.
\bibitem{sheaf} G. Corti\~nas, {\it Periodic cyclic homology as sheaf
cohomology}, $K$-theory {\bf 20} (2000) 175-200.
\bibitem{hanwi} G. Corti\~nas, {\it Cyclic homology of $H$-unital
(pro-)algebras, Lie algebra homology of matrices and a paper of Hanlon's}.
Preprint. Available at http://arxiv.org/abs/math.KT/0504148
\bibitem{val} G. Corti\~nas, C. Valqui, {\it Excision in bivariant
periodic cyclic cohomology: a categorical approach.} $K$-theory
{\bf 30} (2003) 167-201.
\bibitem{cq} J. Cuntz, D. Quillen, {\it Excision in periodic bivariant
cyclic cohomology}, Invent. Math. {\bf 127} (1997), 67-98.
\bibitem{grw} S. Geller, L. Reid, C. Weibel, {\it The cyclic homology
and $K$-theory of curves}, J. reine angew. Math. {\bf 393} (1989) 39-80.
\bibitem{gw0} S. Geller, C. Weibel, {\it $K_1(A,B,I)$}, J. reine angew. Math.
{\bf 342} (1983) 12-34.
\bibitem{gw} S. Geller, C. Weibel, {\it $K(A,B,I)$ II}, $K$-theory {\bf 2}
(1989) 753-760.
\bibitem{gw2} S. Geller, C. Weibel, {\it Hodge decompostion of Loday symbols
in $K$-theory and cyclic homology}, $K$-theory {\bf 8} (1994) 587-632.
\bibitem{gersten} S.M. Gersten, {\it On the spectrum of algebraic $K$-theory},
Bull. Amer. Math. Soc. {\bf 78} (1972) 216-220.
\bibitem{gg} J.A. Guccione, J.J. Guccione, {\it The theorem of excision for Hochschild and cyclic homology.}
J. Pure Appl. Algebra {\bf 106} (1996) 57-60.
\bibitem{pregoo} T. Goodwillie, {\it Cyclic homology, derivations,
and the free loopspace.} Topology {\bf 24} (1985) 187-215.
\bibitem{goo} T. Goodwillie, {\it Relative algebraic $K$-theory and cyclic
homology}, Ann. Math. {\bf 124} (1986), 347-402.
\bibitem{gl} D. Guin-Wal\'ery, J.L. Loday, {\it Obstruction \`a l'excision
en $K$-th\`eorie algebrique}, Springer Lect. Notes Math. {\bf 854} (1981)
179-216.
\bibitem{kav} M. Karoubi, O.E. Villamayor, {\it $K$-th\`eorie
algebrique et $K$-th\`eorie topologique I.} Math. Scand. {\bf 28}
(1971) 265-307.
\bibitem{keu} F. Keune, {\it Doubly relative $K$-theory and the relative
$K_3$}, J. Pure Appl. Algebra {\bf 20} (1981) 39-53.
\bibitem{lod} J. L. Loday, {\it Cyclic homology}, 1st ed.
Grund. math. Wiss. 301. Springer-Verlag Berlin, Heidelberg 1998.
\bibitem{may} J.P. May, {\it Simplicial objects in algebraic topology},
Chicago Lect. Notes in Math, Univ. of Chicago Press, Chicago and London,
1976.
\bibitem{ow} C. Ogle, C. Weibel, {\it Relative $K$-theory and cyclic homology}.
Preprint.
\bibitem{qs} A. Suslin, M. Wodzicki, {\it Excision in algebraic $K$-theory},
Ann. of Math. {\bf 136} (1992) 51-122.
\bibitem{tho} R.W. Thomason, {\it Algebraic $K$-theory and \'etale cohomology}
Ann. scient. \'Ec. Norm. Sup., $4$\'eme  s\'erie, {\bf 18} (1980) 437-552.
\bibitem{chu} C. Weibel, {\it Mayer-Vietoris sequences and mod $p$ $K$-theory},
Springer Lect. Notes Math. {\bf 966} (1982) 390-406.
\bibitem{soval} C. Valqui, {\it Weak equivalence of pro-complexes and excision in topological
Cuntz-Quillen theory}. Preprint SFB 478, Heft 88. Available from

http://wwwmath.uni-muenster.de/math/inst/sfb/about/publ/

\bibitem{wod}
M. Wodzicki, {\it Excision in cyclic homology and in rational algebraic
$K$-theory}, Ann. of Math. {\bf 129} (1989) 591-639
\bibitem{zee}
E.C. Zeeman, {\it A proof of the comparison theorem for spectral sequences,}
Proc. Camb. Philos. Soc. {\bf 53} (1957) 57-62.

\end{thebibliography}
\end{document}